# Numerical Analysis of Automodel Solutions for Superdiffusive Transport

Alexander B. Kukushkin, Vladislav S. Neverov, Petr A. Sdvizhenskii, and Vladimir V. Voloshinov

*Abstract*—The distributed computing analysis of the accuracy of automodel solutions for the Green's function of a wide class of superdiffusive transport of perturbation on a uniform background is carried out. The approximate automodel solutions have been suggested for the 1D transport equation with a model long-tailed step-length probability distribution function (PDF) with various power-law exponents. These PDFs describe the transport dominated by the Lévy flights. Massive computing experiments were done to verify automodel solutions. The Everest distributed computing platform and the cluster at NRC 'Kurchatov Institute' were used. The results verify the high accuracy of automodel solutions in a wide range of space-time variables and suggest extending the developed method of automodel solutions to a wider class of stochastic phenomena.

## I. Introduction

Superdiffusive transport phenomena are of growing interest in the literature and various applications. In the case of normal (or ordinary) diffusion, defined as the Brownian motion described by the differential equation of the Fokker-Planck type, the Green's function is a Gaussian which argument determines the scaling law for the propagation front, $r_{fr}(t) \sim (Dt)^{\beta}$, where $\beta = \frac{1}{2}$ and $D$ is the diffusion coefficient. This law is violated for a broad class of phenomena where the free-path length (step length) provided by the long-tailed, power-law probability distribution function (PDF) leads to the exponent $\beta > \frac{1}{2}$, that is called a superdiffusive transport (see, e.g., [1]-[4]). In this case the dominant contribution to the transport comes from the long-free-path carriers (named, by P. Mandelbrot [5], Lévy flights). In various physics problems, superdiffusion was qualified/named as a non-local transport which is described by an integral, in space variables, equation which is non-reducible to a differential one: for instance, Biberman-Holstein equation [6]-[8] for the space-time evolution of the density of excited atoms/ions produced by the radiative transfer in the spectral lines in gases and plasmas; excitation transport by phonons [9]; (non-stationary) heat transport by the longitudinal, electron Bernstein, waves in plasmas [10]; photon-assisted transport of minority carriers in semiconductors (photo-excited holes in n-type InP) [11].

A wide class of non-stationary superdiffusive transport on a uniform background with a power-law decay, at large distances, of the step-length PDF was shown [12] to possess an approximate automodel solution. The solution for the Green's function was constructed using the scaling laws for the propagation front (relevant-to-superdiffusion average displacement) and asymptotic solutions far beyond and far in advance of the propagation front. These scaling laws were shown to be determined essentially by the long-free-path carriers (Lévy flights). The validity of the suggested automodel solution was proved by its comparison with numerical solutions in the one-dimensional (1D) case of the transport equation with a simple long-tailed PDF with various power-law exponents and in the 3D case of the Biberman-Holstein equation of the resonance radiation transfer for various (Doppler, Lorentz, Voight and Holtsmark) spectral line shapes. The analysis of the limits of applicability of the automodel solution in the 1D case was continued in [13]. However, the full-scale numerical analysis of the limits of applicability was not done yet.

Here, we carry out massive computations of the exact solution of the transport equation with a simple long-tailed PDF with various power-law exponents and compare these results with automodel solutions [12] to identify the limits of their applicability. Preliminary computations have been done on the cluster at NRC 'Kurchatov Institute', http://ckp.nrcki.ru/ by means of the Everest, http://everest.distcomp.org/, a new distributed computing platform supporting the publication, execution and composition of applications running across distributed computing resources [14]. A special Everest application, the so called Parameter Sweep [15], was used to run a bunch of independent tasks at the cluster via special Everest-agent installed on the cluster. However, final computation that includes hundred independent tasks was performed directly at the cluster via SLURM commands.

Manuscript received December 28, 2017. This work was supported in part by the Russian Foundation for Basic Research (projects RFBR-15-07-07850-a, 15-29-07068-ofi-m). This work has been carried out using computing resources of the federal collective usage center Complex for Simulation and Data Processing for Mega-science Facilities at NRC 'Kurchatov Institute', http://ckp.nrcki.ru/

A. B. Kukushkin is with the National Research Centre 'Kurchatov Institute', Moscow, 123182, Russian Federation (e-mail: Kukushkin_AB@nrcki.ru, kukushkin.alexander@gmail.com) and the National Research Nuclear University MEPhI, Moscow, 115409, Russia.
V. S. Neverov, is with the National Research Centre 'Kurchatov Institute', Moscow, 123182, Russia (e-mail: vs-never@hotmail.com).
P. A. Sdvizhenskii is with the National Research Centre 'Kurchatov Institute', Moscow, 123182, Russia (e-mail: sdvinpt@gmail.com).
V. V. Voloshinov is with Institute for Information Transmission Problems (Kharkevich Institute) of Russian Academy of Science, Moscow, 127051, Russia (vv_voloshinov@iitp.ru).

## II. MAIN EQUATIONS

### A. Transport Equation and Exact Solution

We consider the 1D transport on a uniform background, described by the equation for spatial density $f(x, t)$ of an excitation of the background medium, which may evolve due to the exchange of excitation between various points of the medium via emission and absorption of the carriers (here the retardation, caused by the finite velocity of carriers, and the quenching of perturbation are neglected; the derivation of this equation for a possible mechanism of interaction between the medium and the carriers of medium's perturbation is given in Appendix in [12]). The equation for the Green's function $f_G(x, t)$ is as follows:

$$\frac{\partial f(x,t)}{\partial t} = \frac{1}{\tau}\int_{-\infty}^{\infty} W(|x-x_1|)f(x_1,t)\mathrm{d}x_1 - \frac{1}{\tau}f(x,t) + \delta(x)\delta(t) \cdot \quad (1)$$

where $W(x)$ is a step-length PDF (i.e. the probability that the carrier, emitted at some point, is absorbed at a distance $x$ from that point), $1/\tau$ is the absolute value of the emission rate (i.e. $\tau$ is the average waiting time between the absorption and reemission of the carrier), $\delta(x)\delta(t)$ stands for the source function, which is the rate of production of excitation by an external source (i.e. a source which differs from the excitation of the medium due to absorption described by the $W$ function), and, for the Green's function, corresponds to an instant point source. The uniformity of the background assumes that, first, the $W$ is a function of only one variable — the distance between the points of emission and absorption — and, second, $\tau$ is a constant. Hereafter we use the dimensionless time and space coordinate, assuming the normalization of time by $\tau$ and using a dimensionless PDF.

We take the PDF in the following simple form which possesses a long tail and the infinite value of the mean square displacement of the excitation from the source:

$$W(\rho) = \frac{\gamma}{2(1+\rho)^{\gamma+1}}, \quad 0 < \gamma < 2, \quad \rho = |x-x'|,$$
$$\int_{-\infty}^{\infty} W(|x-x'|)dx' = 1. \quad (2)$$

Here, the upper bound for parameter $\gamma$ corresponds to the boundary between the superdiffusive evolution of the Green's function and diffusive, Brownian one.

The analytic solution of (1) with the PDF (2) has the form [12]:

$$f_{exact}(x,t) = \frac{1}{2\pi}\int_{-\infty}^{\infty}\cos(px)\exp\left[-tp\int_{0}^{\infty}\frac{\sin px_1}{(1+x_1)^\gamma}dx_1\right]dp. \quad (3)$$

### B. Automodel Solution and the Algorithm of Verifying the Self-Similarity

The approximate automodel solution of the superdiffusive transport equation (1) for the PDF (2) has the form [12]:

$$f_{auto}(x,t) = t\frac{\gamma}{2\left[1+\rho g\left(\rho_{fr}(t)/\rho\right)\right]^{\gamma+1}}, \quad \rho = |x|, \quad (4)$$

$$\rho_{fr}(t) = (t+1)^{1/\gamma} - 1 \quad (5)$$

where $g$ is a function of a single variable, which has to be reconstructed from comparison with exact solution, while the asymptotic behavior of the function $g$ is known from the scaling laws directly derived from the transport equation (1):

$$g(s) = 1, \ s \ll 1; \quad (6)$$
$$g(s) = \alpha s, \ s \gg 1, \quad \alpha(\gamma) = 2^{1/\gamma}\left[\frac{\gamma\pi}{2}[I(\gamma)]^{1/\gamma}\right]^{1/(\gamma+1)} \quad (7)$$

The procedure of the reconstruction of the function $g$ via comparing the automodel solution (4) with the exact solution (3) is as follows:

$$Q_W(x,t) = \frac{1}{\rho}\left[\left(\frac{\gamma t}{2f_{exact}(x,t)}\right)^{\frac{1}{\gamma+1}} - 1\right], \quad \rho \equiv |x|, \quad (8)$$

$$Q_W(\rho, t(\rho,s)) \equiv Q_W(s,\rho) = g(s), \quad (9)$$
$$Q_W(\rho(t,s),t) \equiv Q_W(s,t) = g(s), \quad (10)$$

where the functions $t(\rho, s)$ and $\rho(t, s)$ are determined by the relation

$$s = ((t+1)^{1/\gamma} - 1)/\rho. \quad (11)$$

Note, that the propagation front (5) is a partial case of the relation defined by (11): $\rho_{fr}(t) = \rho(t, s = 1)$.

## III. DISTRIBUTED COMPUTING PARAMETERS AND IMPLEMENTATION

In order to obtain functions $g(s)$ for different values of $\gamma$, the exact solutions $f_{exact}(\rho,t)$ were calculated. The values of $t$ were evenly spaced on a log scale (100 points per power) in the range from $t_{min} = 30$ to $t_{max} = 10^8$, creating the numerical mesh of total 653 points. The values of $s$ were also evenly spaced on a log scale (50 points per power) in the range from $s_{min} = 0.01$ to $s_{max} = 1000$, creating the numerical mesh of total 501 points. The respective values of $\rho$ are defined by these two numerical meshes, using (11). The values of $\gamma$ were evenly spaced in the range from $\gamma_{min} = 0.5$ to $\gamma_{max} = 1.5$, creating the numerical mesh of total 101 points. We used distributed computing to calculate $f_{exact}(\rho,t)$ on these meshes.

We had to solve 101 independent tasks that corresponds to the mesh of variable $\gamma$. All analytic expressions (3)-(11) have been implemented in the form of a Python program, using the Numpy and SciPy, well-known scientific libraries (www.scipy.org).

The less value of $\gamma$ the more time is spent for computation, e.g., for $\gamma = 0.5$ computations took 8.5 hours on a single core of Intel Xeon E5-2650v2 CPU. Some tests have been performed via the Parameter Sweep, a special Everest application [15], enabling one to run a set of independent tasks at the cluster via special Everest-agent. It was found that the available implementation of the Everest-agent for the SLURM (see agent's sources at https://gitlab.com/everest/agent) has a limit on the number of tasks running in parallel on the cluster (not more than 64). The final computations have been done by the SLURM commands directly from the master host.

The problem was to avoid the default SLURM setting: not more than 64 job per SLURM user (in our case) with one-thread Python process 'inside' each of jobs that run in parallel. The decision was to use the SLURM *srun* command enabling one to run unlimited number of processes within a single SLURM job. The content of appropriate SLURM batch file is shown in Fig. 1. The SLURM-script has been run by the *sbatch* command.

```sh
#!/bin/sh
#SBATCH -D /s/ls4/users/vvvoloshinov/automodel
#SBATCH -n 101
#SBATCH -o %j.out
#SBATCH -e %j.err
#SBATCH -t 9:00:00
#SBATCH -p hpc4-3d

for gamma in `seq 0.5 0.01 1.5`;
 do
  srun --output=$SLURM_JOBID.gamma_${gamma}.out.txt\
       --error=$SLURM_JOBID.gamma_${gamma}.err.txt\
       |-n1 -c1 bash `pwd`/run-selfsim-task.sh $gamma &
 done
wait
```

**Fig. 1.** Example of the SLURM batch file content used to perform the computation for mesh of $\gamma$ from $\gamma_{min} = 0.5$ to $\gamma_{max} = 1.5$ with 0.01 increment. The Python process is called within *run-selfsim-task.sh* Bash-script.

## IV. RESULTS OF NUMERICAL ANALYSIS

Besides direct comparison of the automodel, (4), and exact, (3), solutions of the transport equation (1), an inverse problem is solved to identify the range of variable in the two-dimensional space, {time, space coordinate}, where the automodel solution is accurate within 10% error bars around the exact solution. This boundary, dubbed $t_{10\%}$, as a function of time only, i.e. for the entire range of space coordinate, is shown in Fig. 2.

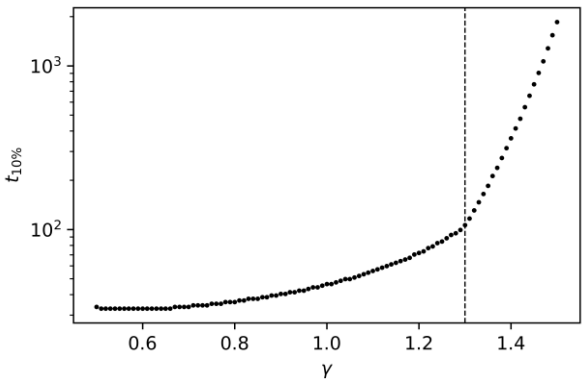

**Fig. 2.** The values of time, $t_{10\%}$, for which the condition $0.9 \leq f_{auto}(\rho,t > t_{min})/f_{exact}(\rho,t > t_{min}) \leq 1.1$ is met for various values of $\gamma$.

The minimal value, $\gamma=0.5$, corresponds to the longest tail of the PDF of the type (2), known in physical problems (see, e.g., [16]). The maximal value, $\gamma=1.5$ corresponds to the case, which is close enough to the boundary between the superdiffusive and diffusive evolution.

Analysis of automodel solutions is illustrated with Figs. 3-5 for three values of $\gamma$. Each case is represented with three figures. The functions $Q_w(s,t)$ (10) are shown for different values of $t$ in the range from $t_{10\%}(\gamma)$ to $t_{max} = 10^8$. The normalized functions $Q_w(s,t)/\{Q_w\}_{av}(s)$, where subscript *av* denotes averaging over time from $t_{min} = 30$ to $t_{max} = 10^8$, and the relative errors of the automodel solution $f_{auto}(\rho,t)/f_{exact}(\rho,t)$ are shown for the same range of time.

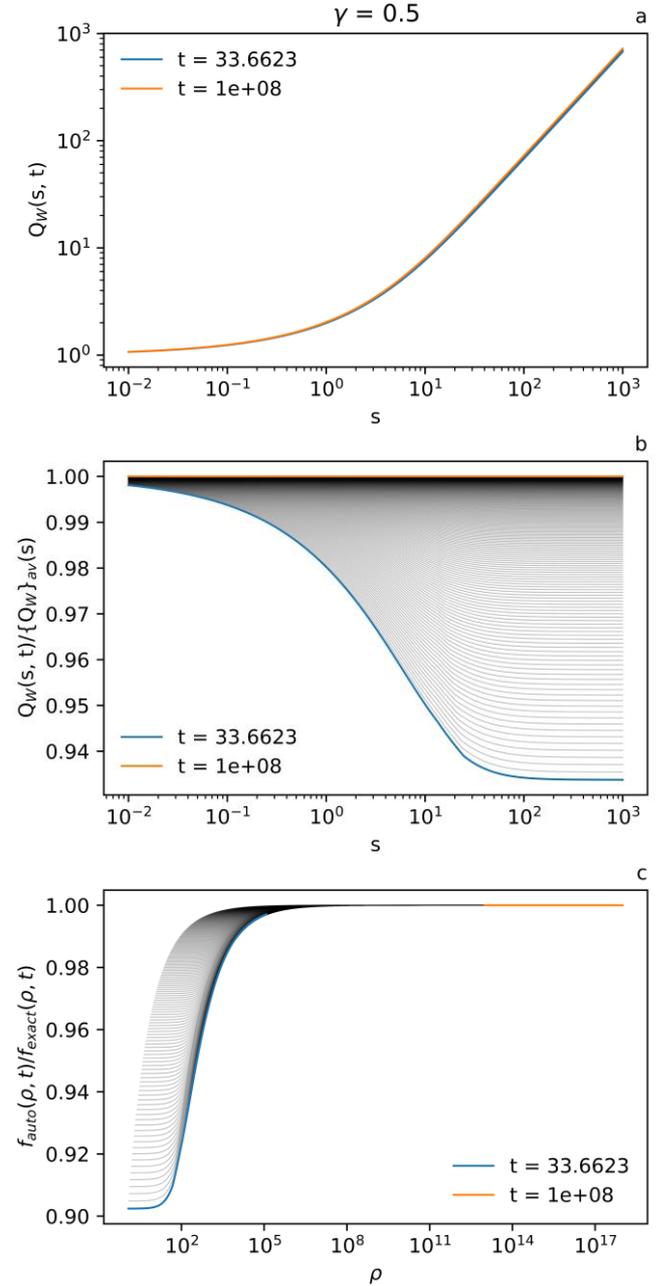

**Fig. 3.** Functions $Q_w(s,t)$ (10) for $\gamma = 0.5$ and different values of $t$ in the range from $t_{10\%} = 33.66$ to $t_{max} = 10^8$ (a). Normalized functions $Q_w(s,t)/\{Q_w\}_{av}(s)$ for the same range of $t$ (b). The relative errors of the automodel solution $f_{auto}(\rho,t)/f_{exact}(\rho,t)$ for the same range of $t$ (c).

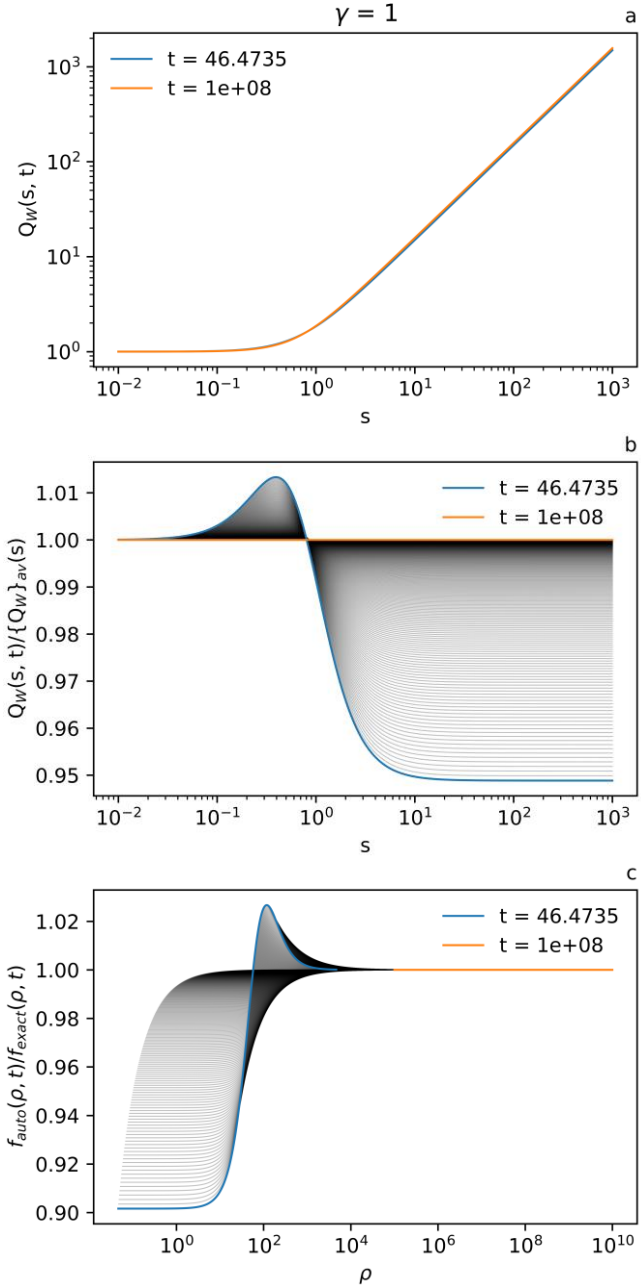

**Fig. 4.** The same plots as in Fig. 3 but for $\gamma = 1.0$ and the values of $t$ in the range from $t_{10\%} = 46.47$ to $t_{max} = 10^8$.

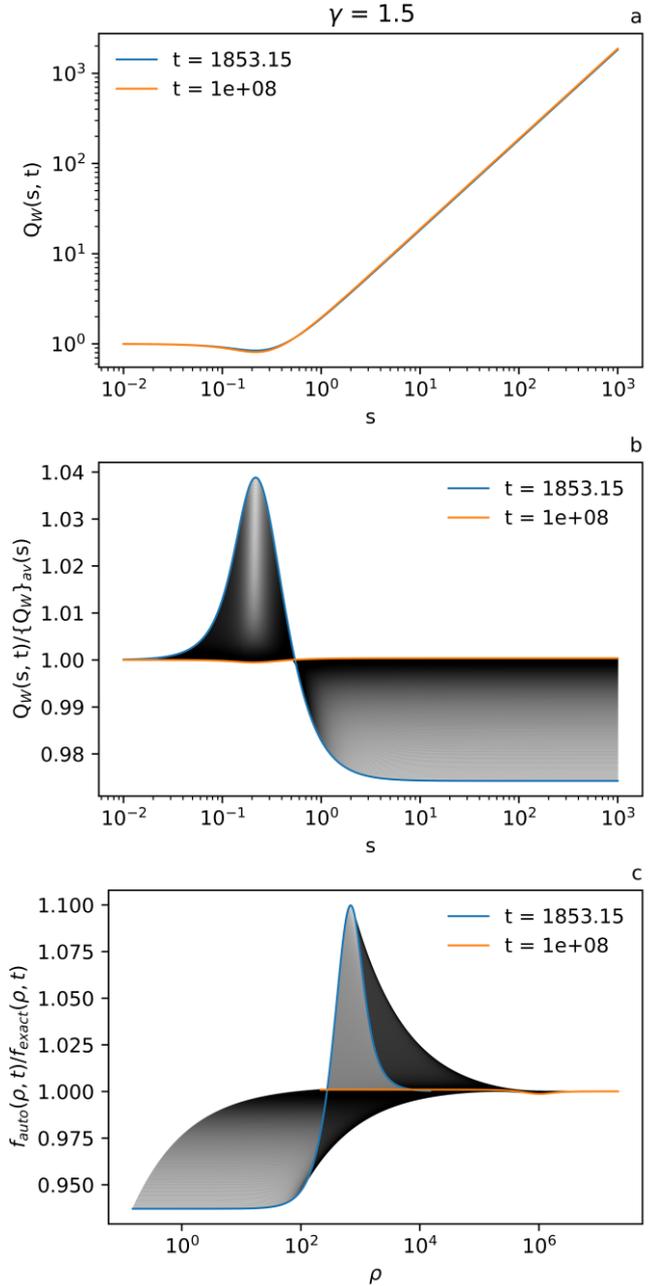

**Fig. 5.** The same plots as in Fig. 3 but for $\gamma = 1.5$ and the values of $t$ in the range from $t_{10\%} = 1853.15$ to $t_{max} = 10^8$.

The results shown in Fig. 3-5 illustrate the following main features of the superdiffusivity as itself and of the automodel solutions.

If the region in $\{t, \rho\}$ space, where 10% is not exceeded, is identified (figures c), the region, where maximal error is located, moves from the point $\{t_{min}, 0\}$ for $\gamma < 1.3$ to the region $s \sim 0.2$ for $\gamma > 1.3$ (figures b). This corresponds to a kink on the curve in Fig. 2, marked with a vertical line. The automodel function $g(s)$, as a function of a single variable, is formed with a high accuracy in the very wide range of $\{t, \rho\}$ space (figures a).

The highest superdiffusivity, being produced by the longest tail of the PDF, needs, as expected, the longest computation time of the exact solution, whereas the applicability of the automodel solutions is limited only by the expected violation at small time. This illustrates the importance of automodel solutions for the most time-consuming problems.

## V. Conclusion

The results, obtained with distributed computing, verified the high accuracy of automodel solutions in a wide range of space-time variables and enabled us to identify the limits of applicability of automodel solutions. All these results suggest extending the developed method of automodel solutions to a wider class of stochastic phenomena.

## Acknowledgment

The authors are grateful to A. P. Afanasiev for the support of collaboration between the NRC "Kurchatov Institute" and